\documentclass[reqno]{amsart}

\usepackage{xcolor}
\usepackage{graphicx}
\usepackage{bbm}
\usepackage{todonotes}
\usepackage{orcidlink}

\usepackage[nobysame,initials,alphabetic]{amsrefs}

\allowdisplaybreaks

\newtheorem{theorem}{Theorem}[section]
\newtheorem{lemma}[theorem]{Lemma}

\numberwithin{equation}{section}

\newcommand{\R}{\mathbb{R}}

\newcommand{\dx}{{\mathrm d}}

\DeclareMathOperator{\erf}{erf}

\title{Central diagonal sections of Gaussian cubes}

\author[F. Fodor]{Ferenc Fodor$^{\orcidlink{0000-0001-9747-1981}}$}

\address{Bolyai Institute, University of Szeged, Aradi v\'ertan\'uk tere 1, 6720 Szeged, Hungary}

\email{fodorf@math.u-szeged.hu}

\author[B. Gonz\'alez Merino]{Bernardo Gonz\'alez Merino$^{\orcidlink{0000-0002-6178-8001}}$}

\address{Departamento de Ingenier\'ia y Tecnolog\'ia de Computadores, \'Area de Matem\'atica Aplicada, Facultad de Inform\'atica, Universidad de Murcia, 30100-Murcia, Spain}

\email{bgmerino@um.es}


\date{\today}

\begin{document}

\begin{abstract}
The investigation of the volume, surface area, and other geometric properties of sections of convex bodies, and in particular cubes, has a long history and a rich literature. However, much less is known when the cube has a volume distribution that is different from the Lebesgue measure; for example, a Gaussian density.
We study the probability densities in the standard cube $B^n_\infty=[-1,1]^n$ of $\R^n$ generated by $e^{-b\|x\|^2}$, $b> 0$. We prove that the limit of the induced Gaussian-type volume of hyperplane sections of $B^n_\infty$ through the origin and orthogonal to a main diagonal is
\[
\sqrt{\frac b\pi}\left (1-4\frac{e^{-b}\sqrt{b}}{2\sqrt{\pi}\erf(\sqrt{b})}\right)^{-\frac12},
\]
as $n\to\infty$. This extends the well-known result of Hensley (1979) for the Lebesgue measure and continues the investigations initiated by Barthe, Gu\'edon, Mendelson, Naor (2005), Zvavitch (2008), and
K\"onig, Koldobski (2013). 
\end{abstract}

\maketitle

\section{Introduction and results}
Let $B^n_\infty=[-1,1]^n$ be the standard $n$-dimensional cube of edge length $2$ centred at the origin. We denote by $\|\cdot\|$ the Euclidean norm in $\R^n$. Let $b\geq 0$ be fixed. For every $x=(x_j)_{j=1}^n\in B_\infty^n$, let 
\[
\dx\gamma_n[b](x)=\frac{e^{-b\|x\|^2}}{\left(\int_{-1}^1 e^{-b\sigma^2}\, \dx \sigma\right)^n}\, \dx x=\frac{\prod_{j=1}^n e^{-bx_j^2}}{\left(\int_{-1}^1 e^{-b\sigma^2}\, \dx \sigma\right)^n}\, \dx x
\]
be the Gaussian-type probability density with parameter $b$ in $B^n_\infty$. Note that if $b=0$, then $\dx \gamma_n[b](x)=1/2^n\dx x$; the uniform density (normalized Lebesgue measure) in $B^n_\infty$.  Let $S^{n-1}$ be the origin-centred unit sphere, and let $\langle\cdot,\cdot\rangle$ denote the Euclidean inner product in $\R^n$. Following K\"onig and Koldobsky \cite{KK13}, we introduce the induced $(n-1)$-measure $\tilde{\gamma}_n[b]$ of the intersection of $B^n_\infty$ with the hyperplane $H(u)=\{x\in \R^n\colon \langle x,u\rangle=0\}$ as follows. For $u\in S^{n-1}$, let
\[
A(u, \gamma_n[b])=
\widetilde{\gamma}_n[b](B^n_\infty\cap H(u)),
\]
where
\[
\widetilde{\gamma}_n[b](B^n_\infty\cap H(u))\colon =\lim_{t\to 0^+} \frac 1{2t}\gamma_n[b](\{x\in B^n_\infty\colon |\langle x, u\rangle|\leq t\}).
\]
K\"onig and Koldobsky \cite{KK13}*{Proposition~2.1} proved \footnote{We note that the product measure $\gamma_n[b]$ in \cite{KK13} misses a factor of $2$ and, as a result, the formula in Proposition 2.1 there also misses a factor of $2$.} that the $\widetilde{\gamma}_n[b]$ measure of hyperplane sections orthogonal to a main diagonal $a=\frac{1}{\sqrt n}(1,\dots,1)$ of $B^n_\infty$ is given by  
\begin{equation}\label{eq:A(a,gamma)}
A(a,\gamma_n[b])
= 
\frac 2\pi \int_0^\infty\left(\frac{\int_0^1\cos(\frac{r}{\sqrt{n}}s)e^{-bs^2}\, \dx s}{\int_0^1e^{-b\sigma ^2}\, \dx \sigma}\right)^n\, \dx r.
\end{equation}

This result generalizes Ball's formula for the Lebesgue measure of central sections of the unit cube \cite{B86}, which corresponds to the special case $b=0$. The origins of Ball's volume formula trace back to P\'olya \cite{P}, see also Bartha, Fodor and Gonz\'alez Merino \cite{BFG21}*{(1)}. 
Using the volume formula for central sections, Ball showed in \cite{B86} that the maximal $(n-1)$-dimensional Lebesgue measure of a hyperplane section of a unit cube is attained precisely when the hyperplane is an $(n-1)$-dimensional subspace that contains an $(n-2)$-dimensional face of $B^n_\infty$; that is, for example, it is parallel to the vector $(1,1,0,\dots,0)$. Ivanov and Tsiutsiurupa \cite{IvTs}, Ambrus \cite{Am}, and Ambrus and G\'argy\'an \cite{AmGa} studied different aspects of local maximizers of central sections of the cube, K\"onig and Rudelson \cite{KoRu} and Moody, Stone, Zach and Zvavitch \cite{MSZZ} investigated non-central sections, and K\"onig and Koldobsky \cite{KK19} dealt with the case of maximizing the surface area. Other aspects of sections of the cube and other convex bodies have recently attracted attention; see, for instance, 
Abel
\cite{Ul18}, 
Alonso--Guti\'errez, Brazitikos and Chasapis
\cite{ABC25}, 
de Loera, Lopez--Campos and Torres
\cite{LLT24}, 
K\"onig
\cites{Kon21, Kon25}, 
K\"onig and Koldobsky
\cite{KonK11},
Lonke
\cite{Lo00}, 
Marichal and Mossinghoff
\cite{MaMo08}, 
Meyer and Pajor
\cite{MePa88}, 
Brandenburg and Meroni
\cite{MeBr25}, 
Nayar and Tkocz
\cite{NaTk23}, 
Pournin
\cites{Po1-23, Po2-23, Po24}.

The problem of finding maximal sections of Gaussian cubes is still open. Barthe, Gu\'edon, Mendelson and Naor \cite{BGMN05} proved general upper bounds for the measure of central hyperplane sections that work for all $b$. Zvavitch \cite{Zvavitch} pointed out that when $b$ is large enough, the central section of the cube, orthogonal to a main diagonal, has a larger $\gamma_n[b]$ measure than the section parallel to the vector $(1,1,0,\dots,0)$. K\"onig and Koldobsky \cite{KK13} quantified Zvavitch's result and proved \cite{KK13}*{Theorem~1.2} that the maximal central sections with respect to the measure $\gamma_n[b]$ are parallel to the vector $(1,1,0,\dots,0)$ if and only if $b<\lambda_0\approx 0.1962627$. Notice that when $b$ is close to $0$, $\gamma_n[b]$ in $B^n_\infty$ is near the Lebesgue measure.

Hensley \cite{H79} proved that the limit of the sequence of the $(n-1)$-dimensional volume of central diagonal sections of $B^n_\infty$ tends to $\sqrt{6/\pi}$ as $n\to\infty$; a result he attributed originally to Selberg. 
K\"onig and Koldobsky \cite{KK19}*{Prop. 6(a)} showed that the volume of central diagonal sections of $B^n_\infty$ is upper bounded by $\sqrt{6/\pi}$. Using Laplace's method and numerical tools, it was established  by Bartha, Fodor and Gonz\'alez Merino \cite{BFG21} that the Lebesgue measure of central sections of $B^n_\infty$, orthogonal to a main diagonal, form a monotonically increasing sequence for $n\geq 3$. We refer to 
Aliev
\cites{Ali20,Ali08}, 
Borwein, Borwein and Leonard
\cite{BBL},
Ron, Ol’hava and Spektor
\cite{KOS} for various properties of the behavior of this sequence. We also note that the volume of central sections can be evaluated explicitly via a closed formula (see Goddard \cite{Go45}, Grimsey \cite{Gr45}, Butler \cite{But60}, Frank and Riede \cite{FR}; see also \cite{BFG21}*{(2)}). For a detailed survey and history on sections of convex bodies, we refer to the paper by Nayar and Tkocz \cite{NaTk23}.

Our main result, Theorem~\ref{thm:main}, is the exact value of the limit of $A(a,\gamma_n[b])$ as $n$ tends to infinity. In particular, \eqref{eq:limit} extends the result of Hensley regarding the volume of central diagonal sections of $B^n_\infty$ mentioned above and can be considered a first step in the investigation of the behavior of the sequence $A(a, \gamma_n[b])$ as $n\to\infty$.

Let $\erf (x)=\frac{2}{\sqrt\pi}\int_0^x e^{-t^2}\, \dx t$ denote the Gaussian error function for $x\in[0,\infty)$.
\begin{theorem}\label{thm:main}
    Let $b>0$. Then
    \begin{equation}
    \lim_{n\rightarrow\infty}A(a,\gamma_n[b]) = 2\sqrt{\frac b\pi}\left (1-4\frac{e^{-b}\sqrt{b}}{2\sqrt{\pi}\erf(\sqrt{b})}\right)^{-\frac12}.\label{eq:limit}
    \end{equation}
\end{theorem}
Notice that the expression of $A(a,\gamma_n[b])$ in \eqref{eq:limit} satisfies $\lim_{n\rightarrow \infty }A(a,\gamma_n[0])=\lim_{b\rightarrow 0^+}\lim_{n\rightarrow\infty}A(a,\gamma_n[b])=\sqrt{\frac6\pi}$, coinciding with the Lebesgue case.

The outline of the proof of Theorem \ref{thm:main} is as follows. We first observe in Section~\ref{sec:LimitCompactIntervals} that over every compact interval $[0,R]$ the sequence of functions in the integrand of \eqref{eq:A(a,gamma)} converges uniformly to an integrable function. Second, in Section~\ref{sec:Tail} we prove that the tails of the sequence of integrals of functions converge to $0$ within intervals $[R,\infty)$ for sufficiently large $R>0$. Combining the two estimates, in Section~\ref{sec:mainproof} we obtain the exact value of the limit of the Gaussian central diagonal sections of the cube $B^n_\infty$.

\section{The limit over compact intervals}\label{sec:LimitCompactIntervals}
	For $b>0$, $n\geq 2$ and $r\in [0,\infty)$, let 
	\[
	f(r,n,b) := \frac{\int_0^1\cos(\frac{r}{\sqrt{n}}s)e^{-bs^2}\, \dx s}{\int_0^1e^{-b\sigma^2}\, \dx \sigma}=\frac{\int_0^1\cos(\frac{r}{\sqrt{n}}s)e^{-bs^2}\, \dx s}{\frac{\sqrt{\pi } \text{erf}\left(\sqrt{b}\right)}{2 \sqrt{b}}}.
	\]

\begin{theorem}\label{lem:UniformConvergenceToLimit}
	For any fixed $b>0$ and $r\in [0,\infty)$, 
	\begin{equation*}
		\lim_{n\to\infty} f^n(r,n,b)
		=\exp\left \{\frac{r^2}{4b} \left(\frac{2 \sqrt{b} e^{-b}}{\sqrt{\pi }
				\erf\left(\sqrt{b}\right)}-1\right)\right\},
	\end{equation*}
    and convergence is uniform in any interval $[0,r_0]$ for any $r_0>0$.
\end{theorem}

\begin{proof}
	 Consider the second order Taylor series of $f(r,n,b)$ with respect to $r$ centred at $r=0$, with the Lagrange error term:
	\[
	f(r,n,b)=1+\frac{r^2 \left(\frac{2 \sqrt{b} e^{-b}}{\sqrt{\pi }
			\erf\left(\sqrt{b}\right)}-1\right)}{4 b n}
		+\frac{\partial^3 f}{\partial r^3}(\xi_{r,n},n,b)\frac{r^3}{3!},
	\]
	for some $\xi_{r,n}\in [0, r]$. Observe that $\xi_{r,n}$ also depends on $b$, but since $b$ is fixed in this argument, we suppress this dependence. By direct calculation, we obtain 
	\begin{align*}
	\frac{\partial^3 f}{\partial r^3}(r,n,b)&
		=\left (\frac{e^{-b} \left(4 \sqrt{b} \sqrt{n} \left(r^2-4 b (b+1) n\right) \sin \left(\frac{r}{\sqrt{n}}\right)-8 b^{3/2} n r \cos \left(\frac{r}{\sqrt{n}}\right)\right)}{\sqrt{\pi }16 b^3 \erf\left(\sqrt{b}\right)n^{\frac 32}}\right.\\
		&\qquad-\left.\frac{r  \left(r^2-6 b n\right)f(r,n,b)}{8 b^3 n^{\frac 32}}\right )\frac{1}{n^{\frac 32}}\\
		&=(h_1(r,n,b)-h_2(r,n,b))\frac{1}{n^{\frac 32}},
	\end{align*}
	where
	\begin{align*}
		h_1(r,n,b)&=\frac{e^{-b} \left(4 \sqrt{b} \sqrt{n} \left(r^2-4 b (b+1) n\right) \sin \left(\frac{r}{\sqrt{n}}\right)-8 b^{3/2} n r \cos \left(\frac{r}{\sqrt{n}}\right)\right)}{\sqrt{\pi }16 b^3 \erf\left(\sqrt{b}\right)n^{\frac 32}},\\
		h_2(r,n,b)&=\frac{r  \left(r^2-6 b n\right)f(r,n,b)}{8 b^3 n^{\frac 32}}.
	\end{align*}
	Thus, by the triangle inequality.
	\begin{align*}
	\left |\frac{\partial^3 f}{\partial r^3}(r,n,b)\right |=\left |(h_1(r,n,b)-h_2(r,n,b))\right | \frac{1}{n^{\frac 32}}
	\leq \left (\left |(h_1(r,n,b)\right |+\left |h_2(r,n,b))\right |\right ) \frac{1}{n^{\frac 32}}.	
	\end{align*}
	
	Let $r_0\in [0,\infty)$ be arbitrary but fixed. If $b>0$, $n\geq 2$ are fixed, then 
	both $h_1$ and $h_2$ are bounded on $[0,r_0]$ as they are continuous functions. Furthermore, if $b>0$ is fixed, then
	\begin{align*}
		|h_1(r,n,b)|&\leq C_1(r_0,b), \quad \text{ and }\\
		|h_2(r,n,b)|&\leq C_2(r_0,b)
	\end{align*}
	for all $n\geq 2$ and $r\in [0,r_0]$ for some suitable $C_1(r_0,b)$ and $C_2(r_0,b)$, meaning that they are uniformly bounded on $[0,r_0]$. The first of these inequalities is clear as the $\sin x$ and $\cos x$ functions are bounded. For the second inequality, note that $f$ is continuous on $[0,r_0]$ for any $n\geq 2$, and thus bounded. Moreover, for any fixed $b>0$, the functions $f(r,n,b)$ are clearly bounded for all $n\geq 2$ uniformly on $[0,\infty)$. Therefore, $h_2$ are $O(n^{-\frac 12})$ as $n\to\infty$ on $[0,r_0]$.
	
	In summary, it follows that for any fixed $b>0$, $r_0\in [0,\infty)$ and $r\in [0,r_0]$, 
		\[
		\left |\frac{\partial^3 f}{\partial r^3}(\xi_{r,n},n,b)\right |\leq (C_1(r_0,b)+C_2(r_0,b))\frac{1}{n^{\frac 32}}\\
		=C(r_0,b)\frac{1}{n^{\frac 32}}
		\]
		for all $n\geq 2$.
		
		Therefore, for any fixed $b$, $r_0\in [0,\infty]$, and $r\in [0,r_0]$, we obtain
		\begin{align*}
		\lim_{n\to\infty} f^n(r,n,b)&=\lim_{n\to\infty}\left (1+\frac{r^2 \left(\frac{2 \sqrt{b} e^{-b}}{\sqrt{\pi }
				\erf\left(\sqrt{b}\right)}-1\right)}{4 b n}
			+\frac{\partial^3 f}{\partial r^3}(\xi_{r,n},n,b)\frac{r^3}{3!}\right )^n\\
		&=\lim_{n\to\infty}\left (1+\frac{r^2 \left(\frac{2 \sqrt{b} e^{-b}}{\sqrt{\pi }
				\erf\left(\sqrt{b}\right)}-1\right)}{4 b n}
		+C(r_0,b)\frac{r_0^3}{n^{\frac 32}3!}\right )^n\\	
		&=\exp\left \{\frac{r^2 \left(\frac{2 \sqrt{b} e^{-b}}{\sqrt{\pi }
				\erf\left(\sqrt{b}\right)}-1\right)}{4 b }\right \}.
		\end{align*}
Now, we show that the convergence is uniform in any  interval $[0,r_0]$.
Let 
\[
a_n(u)=\left (1+\frac un+\frac{C}{n^{\frac 32}}\right )
\]
for some constant $C$. Then
\[
\log a_n(u)=n\log \left (1+\frac un+\frac{C}{n^{\frac 32}}\right ).
\]
Let $x_n(u)=\frac un+\frac{C}{n^{\frac 32}}$. Using the Taylor series of $\log (1+x)$ at $x=0$, we obtain
\[
\log a_n(u)=n\left (x_n(u)+\frac{x_n^2(u)}{2}+O\left (x_n^2(u)\right )\right )=u+\frac{C}{\sqrt n}+O\left (\frac 1n\right ).
\]
Now, by exponentiation and using the Taylor series of $e^x$ at $x=0$, we obtain 
\[
a_n(u)=e^u e^{\frac{C}{\sqrt n}+O\left (\frac 1n\right )}=e^u \left (1+\frac{C}{\sqrt n}+O\left (\frac 1n\right )\right ).
\]
Thus,
\[
\left |a_n(u)-e^u\right |=e^u \left |\frac{C}{\sqrt n}+O\left (\frac 1n\right )\right |\leq e^{u_0} \left |\frac{C}{\sqrt n}+O\left (\frac 1n\right )\right |, 
\]
from which the uniform convergence on $[0,u_0]$ follows. 
\end{proof}		
		

\section{The tail estimate}\label{sec:Tail}
\begin{theorem}\label{thm:uniformStability}
    Let $b>0$, $t\geq 0$.
    Then for every $\varepsilon>0$, there exist $R(\varepsilon)>0$ and $n(\varepsilon)\in\mathbb N$ such that
    \[
    \int_{R(\varepsilon)}^\infty |f(r,n,b)^n|\, \dx r \leq \varepsilon
    \]
    for every $n\geq n(\varepsilon)$.
\end{theorem}

For every $b>0$ and $t\geq 0$ let 
\[
G_b(t):=f(t,1,b)=\frac{\int_0^1\cos(ts)e^{-bs^2}\, \dx s}{\int_0^1e^{-b\sigma ^2}\, \dx \sigma}=\int_0^1\cos(ts)\, \dx \mu_b(s) =\frac{F_b(t)}{N_b},
\]
where
\[
\dx\mu_b(s):=
\frac{e^{-bs^2}\, \dx s}{\int_0^1e^{-b\sigma^2}\, \dx \sigma} \mathbf{1}_{[0,1]}(s),\quad F_b(t):=\int_0^1\cos(ts)e^{-bs^2}\, \dx s,
\]
and
\[
N_b:=\int_0^1e^{-b\sigma^2}\, \dx \sigma.
\]

Notice that $G_b(t)$ is a multiple of the real part of the Fourier transform of the Gaussian $e^{-bs^2}$. The Riemann--Lebesgue Lemma ensures $G_b(t)\leq O(\frac1t)$. However, the actual constant also depends on $b$, and we need a more precise estimate for our arguments.
\begin{lemma}\label{lem:larget}
    Let $b>0$ and $t\geq 0$. Then
    \[
    |G_b(t)| \leq 
    \begin{cases}
        \frac{e}{t},  &  \text{ if }b\in(0,1),\\
       1.34\frac{\sqrt{b}}{t}, & \text{ if }b\geq 1.
    \end{cases}
    \]
\end{lemma}

\begin{proof}
    We start bounding $F_b(t)$ from above. We integrate by parts with $u=e^{-bs^2}$ and $\dx v=\cos(ts)\dx s$. Thus, $\dx u=-2bse^{-bs^2}\dx s$, $v=\frac{\sin(ts)}{t}$ and
    \begin{align*}
        \left | F_b(t)\right |&= \left | \frac 1t \left (e^{-bs^2}\sin (ts) \right)_0^1+\frac{2b}{t}\int_0^1 s e^{-bs^2}\sin (ts)\, \dx s\right |\\
        &\leq \left | \frac 1t e^{-b}\sin (t) \right |+\frac{2b}{t}\left |\int_0^1 s e^{-bs^2}\sin (ts)\, \dx s\right |\\
        &\leq\frac 1t e^{-b}+\frac{2b}{t}\int_0^1 s e^{-bs^2}\, \dx s\\
        &= \frac 1t e^{-b}-\frac{2b}{t}\frac{e^{-b}-1}{2b}\, \dx s=\frac 1t.
    \end{align*}

Next, we bound $N_b$ from below. Notice that $N_b$ is a strictly monotonically decreasing function of $b$.

If $b\geq 1$, then using the substitution $u=\sqrt b \sigma$, $\dx u=\sqrt b\, \dx \sigma$, we get
\begin{align*}
    N_b=\frac{1}{\sqrt b}\int_0^{\sqrt b} e^{-u^2}\, \dx u\geq \frac{1}{\sqrt b}\int_0^{1} e^{-u^2}\, \dx u=\frac{1}{\sqrt b}N_1\geq \frac{0.7468}{\sqrt b}, 
\end{align*}
where we used $b\geq 1$. Thus,
\[
N_b\geq 
\begin{cases}
    0.7468>e^{-1}, & \text{ if } b\in (0,1),\\
    \frac{0.7468}{\sqrt b}, & \text{ if } b\geq 1.
\end{cases}
\]
Since $G_b(t)=F_b(t)/N_b$, the statement of the lemma follows.

\end{proof}

\begin{lemma}\label{lem:smallt}
Let $t\geq 0$. Then
\begin{equation*}
|G_b(t)|\leq
\begin{cases}
    1-\frac{1}{16}t^2, & \text{ if } b\in (0,1) \text{ and } t\leq 2.75\\
    1-\frac{1}{25b}t^2, & \text{ if } b\geq 1 \text{ and } t\leq 1.4\sqrt b.
\end{cases}
\end{equation*}
\end{lemma}

\begin{proof}
    Notice that $\cos x\leq 1-\frac{x^2}{4}$ for $x\in [0,2.75]$. If $s\in (0,1)$ and $t\leq 2.75$, or $b\geq 1$, $t\leq 1.4\sqrt b$ and $s\leq 1/\sqrt b$, then
    \[
    \cos (ts)\leq 1-\frac{t^2s^2}{4}.
    \]
Now we split the domain of integration in $G_b(t)$ at $1/\sqrt{b}$ (note that if $b\geq 1$ then the second term is void as $\mu_b(s)$ is supported on $[0,1]$)
\begin{align*}
    |G_b(t)|&=\left |\int_0^{\frac{1}{\sqrt b}} \cos(ts)\, \dx \mu_b(s)+\int_{\frac{1}{\sqrt b}}^1 \cos(ts)\, \dx \mu_b(s) \right |\\
    &\leq \left |\int_0^{\frac{1}{\sqrt b}} 1-\frac{t^2s^2}{4}\, \dx \mu_b(s) \right |+\int_{\frac{1}{\sqrt b}}^1 \left |\cos(ts)\right |\, \dx \mu_b(s)\\
    &\leq 1-\frac{t^2}{4}\int_0^{\frac{1}{\sqrt b}} s^2\, \dx \mu_b(s).
\end{align*}

If $b\in (0,1)$, then
\[
|G_b(t)|\leq 1-\frac{t^2}{4}\int_0^1 s^2\, \dx \mu_b(s)
\]
We now use that the moment of second order 
\[
\int_0^1 s^2\, \dx \mu_b(s) = \frac{\int_0^1s^2e^{-bs^2}\,\dx s}{\int_0^1e^{-b\sigma^2}\,\dx \sigma}
\]
of $e^{-bs^2}$ is decreasing in $b$. Indeed, let
\[
B(b):=\int_0^1s^2e^{-bs^2}\, \dx s \quad\text{and}\quad R(b):=\frac{B(b)}{N_b}.
\]
Differentiating $R(b)$ with respect to $b$ we obtain
\[
R'(b) = \frac{B(b)^2-N_bC(b)}{N_b^2},
\]
where $C(b):=\int_0^1s^4e^{-bs^2}\, \dx s$. The Cauchy-Schwarz inequality yields 
\[
C(b)^2 = \left(\int_0^1e^{-\frac{bx^2}{2}}\cdot (x^2e^{-\frac{bx^2}{2}})\, \dx x\right)^2 
\leq \int_0^1 e^{-bx^2}\, \dx x \int_0^1 x^4e^{-bx^2}\, \dx x = N_bC(b).
\]
Hence $R'(b)<0$ and therefore $R(b)$ is decreasing in $b$.
Thus, we obtain
\[
|G_b(t)| \leq 1-\frac{t^2}{4} \int_0^1 s^2\, \dx \mu_1(s) \leq 1-\frac{0.25}{4}t^2.
\]

If $b\geq 1$, then $1/\sqrt{b}\leq 1$ and we consider
    \[
    \int_0^\frac{1}{\sqrt{b}} s^2\, \dx \mu_b(s) = \frac{\int_0^\frac{1}{\sqrt{b}}s^2e^{-bs^2}\, \dx s}{\int_0^1e^{-bs^2}\, \dx s}.
    \]
We bound the numerator from below and the denominator from above. Regarding the former, since $s\leq1/\sqrt{b}$, then $e^{-bs^2} \geq e^{-1}$, and thus
    \[
    \int_0^\frac{1}{\sqrt{b}}s^2e^{-bs^2}\, \dx s \geq e^{-1}\int_0^\frac{1}{\sqrt{b}}s^2\, \dx s = \frac{e^{-1}}{3b^{\frac{3}{2}}}.
    \]
    For the denominator, we obtain
    \[
    \int_0^1e^{-b\sigma^2}\, \dx \sigma \leq \int_0^\infty e^{-b\sigma^2}\, \dx \sigma = \frac{\sqrt{\pi}}{2\sqrt{b}}.
    \]
 Therefore,
    \[
    \int_0^\frac{1}{\sqrt{b}} s^2\, \dx \mu_b(s) \geq \frac{2e^{-1}}{3\sqrt{\pi}b},
    \]
    and hence
    \[
    \begin{split}
    |G_b(t)| 
    & \leq 1-\frac{t^2}{4}\frac{2e^{-1}}{3\sqrt{\pi}b} \leq 1-\frac{1}{25}\frac{t^2}{b}.
    \end{split}
    \]

\end{proof}

\begin{proof}[Proof of Theorem \ref{thm:uniformStability}]
    We start with proving the case $b\in(0,1)$. Let $R>0$ and assume that $2.73\sqrt n\geq R$. We split the expression \eqref{eq:A(a,gamma)} into
    \begin{equation}\label{eq:b01Integral}
    \int_R^\infty \left|f(r,n,b)^n\right|\, \dx r = \int_R^{2.73\sqrt{n}} \left|f(r,n,b)^n\right|\, \dx r + \int_{2.73\sqrt{n}}^\infty \left|f(r,n,b)^n\right|\, \dx r.
    \end{equation}
    We start with bounding the first term. Since $t=\frac{r}{\sqrt{n}}\leq 2.73$, by Lemma \ref{lem:smallt} we have that
    \[
    \begin{split}
    \int_R^{2.73\sqrt{n}} \left|f(r,n,b)^n\right|\, \dx r & \leq \int_R^{2.73\sqrt{n}} \left(1-\frac{1}{16}\frac{r^2}{n}\right)^n\, \dx r \\
    & \leq \int_R^{2.73\sqrt{n}} e^{-\frac{1}{16}r^2}\, \dx r.
    \end{split}
    \]
    Evidently, for every $\varepsilon>0$, we can choose $R(\varepsilon)>0$ large enough so that 
    \[
    \int_{R(\varepsilon)}^{\infty} e^{-\frac{1}{16}r^2}\, \dx r \leq \frac{\varepsilon}{2}.
    \]
    We now bound the second term. Since $t=\frac{r}{\sqrt{n}} \geq 0$, Lemma \ref{lem:larget} yields
    \[
    \int_{2.73\sqrt{n}}^\infty \left|f(r,n,b)^n\right|\, \dx r 
    \leq \int_{2.73\sqrt{n}}^\infty \left(\frac{e}{r/\sqrt{n}}\right)^n\, \dx r 
    = \frac{2.73\sqrt{n}}{n-1}\left(\frac{e}{2.73}\right)^n.
    \]
    Again, since $e/2.73<1$ there exists $n(\varepsilon)\in\mathbb N$ such that 
    \[
    \frac{2.73\sqrt{n}}{n-1}\left(\frac{e}{2.73}\right)^n\leq \frac{\varepsilon}{2}
    \]
    for every $n\geq n(\varepsilon)$. Joining the left and the right estimates tells us that \eqref{eq:b01Integral} is upper bounded by $\varepsilon$ when choosing $R:=R(\varepsilon)$ and for every $n\geq n(\varepsilon)$ with $n(\varepsilon)\geq \frac{R(\varepsilon)^2}{2.73^2}$, concluding the case $b\in(0,1)$.

We now finish with the case $b \geq 1$. Let $R>0$ be the fixed lower value of the integral domain and assume that $1.37\sqrt b\sqrt n\geq R$. We then split the expression \eqref{eq:A(a,gamma)} onto
    \begin{equation}\label{eq:b1InftyIntegral}
    \int_R^\infty \left|f(r,n,b)^n\right|\, \dx r = \int_R^{1.37\sqrt{b}\sqrt{n}} \left|f(r,n,b)^n\right|\, \dx r + \int_{1.37\sqrt{b}\sqrt{n}}^\infty \left|f(r,n,b)^n\right|\, \dx r.
    \end{equation}
    We start bounding the left hand side. Since $t=\frac{r}{\sqrt{n}}\leq 1.37\sqrt{b} < 1.4\sqrt{b}$, by Lemma \ref{lem:smallt} we have that
    \[
    \begin{split}
    \int_R^{1.37\sqrt{b}\sqrt{n}} \left|f(r,n,b)^n\right|\, \dx r & \leq \int_R^{1.37\sqrt{b}\sqrt{n}} \left(1-\frac{1}{25}\frac{r^2}{bn}\right)^n\, \dx r \\
    & \leq \int_R^{1.37\sqrt{b}\sqrt{n}} e^{-\frac{1}{25b}r^2}\, \dx r.
    \end{split}
    \]
    Evidently, for every $\varepsilon>0$ and fixed $b\geq 1$, we can choose $R(\varepsilon)>0$ large enough so that 
    \[
    \int_{R(\varepsilon)}^{\infty} e^{-\frac{1}{25b}r^2}\, \dx r \leq \frac{\varepsilon}{2}.
    \]
    We now proceed with the second term. Notice that $t=\frac{r}{\sqrt{n}} \geq 1.37\sqrt{b} \geq\sqrt{b}$ and thus Lemma \ref{lem:larget} we get that
    \[
    \int_{1.37\sqrt{b}\sqrt{n}}^\infty \left|f(r,n,b)^n\right|\, \dx r 
    \leq \int_{1.37\sqrt{b}\sqrt{n}}^\infty (1.34\sqrt{b})^n\frac{\sqrt{n}^n}{r^n}\, \dx r 
    = \frac{1.37\sqrt{b}\sqrt{n}}{n-1}\left(\frac{1.34}{1.37}\right)^n.
    \]
    Again, since $1.34/1.37<1$ and $b\geq 1$ is fixed, there exists $n(\varepsilon)\in\mathbb N$ such that 
    \[
    \frac{1.37\sqrt{b}\sqrt{n}}{n-1}\left(\frac{1.34}{1.37}\right)^n\leq \frac{\varepsilon}{2}
    \]
    for every $n\geq n(\varepsilon)$. Joining the left and the right estimates tells us that \eqref{eq:b1InftyIntegral} is upper bounded by $\varepsilon$ when choosing $R:=R(\varepsilon)$ and for every $n\geq n(\varepsilon)$ with $n(\varepsilon)\geq \frac{R(\varepsilon)^2}{1.37^2b}$, concluding the case $b\geq 1$ and the theorem.
\end{proof}

\section{Proof of Theorem~\ref{thm:main}}\label{sec:mainproof}
We now show Theorem \ref{thm:main}. Let us define
\[
g_b(r):=\exp\left \{\frac{r^2}{4b} \left(\frac{2 \sqrt{b} e^{-b}}{\sqrt{\pi }
				\text{erf}\left(\sqrt{b}\right)}-1\right)\right\}.
\]
For a suitable $R>0$, let us write
\[
\begin{split}
& \left|\int_0^\infty f^n(r,n,b)\, \dx r - \int_0^\infty g_b(r)\, \dx r\right| \\ 
& = \left|\int_0^R f^n(r,n,b)\, \dx r + \int_R^\infty f^n(r,n,b)\, \dx r - \int_0^R g_b(r)\, \dx r-\int_R^\infty g_b(r)\, \dx r \right| \\
& \leq \int_0^R \left|f^n(r,n,b)- g_b(r)\right|\, \dx r + \int_R^\infty \left|f^n(r,n,b)\right|\, \dx r+\left|\int_R^\infty g_b(r)\, \dx r\right|.
\end{split}
\]
First, notice that by Theorem \ref{thm:uniformStability}, there exists $R_1(\varepsilon)>0$ and $n_1(\varepsilon)>0$ such that the second term is not larger than $\frac{\varepsilon}{3}$ when choosing $R:=R_1(\varepsilon)$ and every $n\geq n_1(\varepsilon)$. 
Second, 
it is clear that
\[
\int_R^\infty g_b(r) dr= \int_R^\infty \exp\left \{\frac{r^2}{4b} \left(\frac{2 \sqrt{b} e^{-b}}{\sqrt{\pi }
				\text{erf}\left(\sqrt{b}\right)}-1\right)\right\}\, \dx r \leq \frac{\varepsilon}{3}
\]
for some large enough $R=R_2(\varepsilon)$. Take $R:=\max\{R_1(\varepsilon),R_2(\varepsilon)\}$, so that the second and third terms are not larger than $\frac{\varepsilon}{3}$ each whenever $n\geq n_1(\varepsilon)$ (notice that the third term do not depend on $n$ anyways). 
Third, 
notice that Theorem \ref{lem:UniformConvergenceToLimit} ensures the uniform convergence of $f^n(r,n,b)$ to $g_b(r)$ on $r\in[0,R]$, for every $R>0$ when $n\rightarrow\infty$. Thus, for our choice of $R$, there exists $n_2(\varepsilon)\in\mathbb N$ such that for every $n\geq n_2(\varepsilon)$ then the first term in the equation above is not larger than $\frac{\varepsilon}{3}$. Letting $n_0:=\max\{n_1(\varepsilon),n_2(\varepsilon)\}$, we conclude that for every $n\geq n_0$ then 
\[
\left|\int_0^\infty f^n(r,n,b)\, \dx r - \int_0^\infty g_b(r)\, \dx r\right|\leq\varepsilon,
\]
concluding the proof.

\section*{Funding sources}
F. Fodor's research was supported by NKFIH project no.~150151. Project no.~150151 has been implemented with the support provided by the Ministry of Culture and Innovation of Hungary from the National Research, Development and Innovation Fund, financed under the ADVANCED\_24 funding scheme.

B.~Gonz\'alez Merino was partially supported by Ministerio de Ciencia, Innovación y Universidades project PID2022-136320NB-I00/AEI/10.13039/501100011033/FEDER, UE.

\end{document}